\documentclass[11pt]{article}
\usepackage{amsfonts,amssymb,amsmath}\usepackage{verbatim}
\usepackage{graphicx}%\input{jqdefs}
\input{xypic}
\usepackage{longtable}
\usepackage{lscape}
\def\Z{\mathbb Z}\def\Q{\mathbb Q}\def\Qb{\overline{\mathbb Q}}\def\GL{\operatorname{GL}}
\def\End{\operatorname{End}}\def\Gal{\operatorname{Gal}}\def\Frob{\operatorname{Frob}}\def\Br{\operatorname{Br}}
\def\Res{\operatorname{Res}}\def\Disc{\operatorname{Disc}}
\title{Package description and tables for the paper\\
    Fields of definition of building blocks}
\author{Jordi Quer}
\date{May, 31 2006}

\begin{document}\maketitle

\begin{abstract}
This document contains additional information for
the paper: J. Quer. {\sl Fields of definition of building blocks},
Math. Comp. {\bf 78} (2009), no. 265, pp. 537--554.

It contains a description of the implementation in {\tt Magma}
of the algorithms and computational techniques described in that paper,
and also several tables with data that were elaborated using this implementation.
\end{abstract}
\tableofcontents

%---------------------------------------------------------------------------------------------------
\section{Background: definitions and notation}
%---------------------------------------------------------------------------------------------------

Let $f=\sum_{n=1}^\infty a_n q^n\in S_2^{\operatorname{new}}(N,\varepsilon)$
be a newform of weight $2$, level $N$ and Nebentypus $\varepsilon$. Let
$$E=\Q(\{a_n\}_{n\geq 1})\qquad\text{and}\qquad F=\Q(\{a_p^2/\varepsilon(p)\}_{p\nmid N})$$
be the number fields generated, respectively, by all the Fourier coefficients of $f$
and by the numbers
$$\mu_p=a_p^2/\varepsilon(p),$$
for all primes $p$ not dividing the level $N$.

The abelian variety $A_f/\Q$ attached by Shimura to $f$ is a variety of $\GL_2$-type
with $\End^0(A_f)\simeq E$; the isomorphism is given by the identification of
(the restriction to $A_f$ of) the Hecke operator $T_p$ acting on the jacobian $J_1(N)$
with the algebraic integer $a_p$.
The variety $A_f$ factors up to isogeny over $\Qb$ as a power $B_f^r$ of a simple variety $B_f$.

The variety $A_f$ has complex multiplication if, and only if, $B_f$ is a CM elliptic curve.
This happens exactly when the form $f$ has complex multiplication by a quadratic field,
meaning that there exists a nontrivial primitive Dirichlet character $\chi$
(which is necessarily unique, quadratic and odd) such that
$$a_p=\chi(p)\,a_p\qquad\text{for all}\quad p\nmid N.$$
In that case, the curve $B_f$ is
just the curve with complex multiplication by (an order of) that quadratic field
(which is unique up to isogeny).

\def\D{{\mathcal D}}
From now on we assume that $f$ has no complex multiplication;
in this case the varieties $B_f$ are called \emph{building blocks}.
Now, the field $F$ is totally real and the field $E$ is an abelian extension of it
that is either a totally real field when $\varepsilon$ is trivial
(which is equivalent to $A_f$ being a factor of $J_0(N)$ up to $\Q$-isogeny),
or a CM field otherwise.
The endomorphism algebra $\D=\End^0(B_f)$ is a division algebra;
it can be either equal to the totally real field $F$, with $\dim(B_f)=[F:\Q]$
(and $r=[E:F]$ in this case) or to a quaternion algebra over the totally real field $F$
(which is necessarily totally indefinite), with $\dim(B_f)=2[F:\Q]$
(and $r=\frac12[E:F]$ in this case).

For every element $s\in G_F$ there exists a unique primitive Dirichlet character $\chi_s$,
which only depends on the action of $s$ on the field $E$, such that
\begin{equation}\label{definition-phi}
{}^s a_p =\chi_s(p)\, a_p\qquad\text{for all}\quad p\nmid N.
\end{equation}
The character $\chi_s$, and sometimes also the pair $(s,\chi_s)$,
is called an \emph{inner twist} of the form $f$.

Let $F_\delta=F(\{\sqrt{\mu_p}\}_{p\nmid N})$ be the extension of $F$ generated by the
square roots of the elements $\mu_p\in F$, which is a polyquadratic extension of $F$.
For every element $s\in G_F$ there exists a unique primitive quadratic Dirichlet character
$\psi_s$,
which only depends on the action of $s$ on the field $F_\delta$, such that
\begin{equation}\label{definition-psi}
{}^s \sqrt{\mu_p}=\psi_s(p)\,\sqrt{\mu_p}\qquad\text{for all}\quad p\nmid N.
\end{equation}
The character $\psi_s$, and also the pair $(s,\psi_s)$,
will be called a \emph{quadratic degree character} of the form $f$.
The group $\Psi$ of all the quadratic degree characters is an abelian group of exponent $2$,
isomorphic to $\Gal(F_\delta/F)$.
It depends only on the $\Qb$-isogeny class of the building block $B_f$.

The inner twists and quadratic degree characters are related by the following identity
\begin{equation}\label{relation-chi-psi}
\chi_s(p)=\psi_s(p)\sqrt{\varepsilon(p)}^{s-1}\qquad\text{for all}\quad p\nmid N.
\end{equation}
We observe that since the $\chi_s$ depend only on $s|_E$
and the $\psi_s$ depend on $s|_{F_\delta}$,
everything in this identity depends only on the restriction of $s$ to the field
$E\cdot F_\delta=E(\sqrt{\varepsilon})$, where $E(\sqrt{\varepsilon})$ denotes the
quadratic field obtained by adjoining to $E$ the square roots of the values of the
character $\varepsilon$, which is either equal to $E$ or to a quadratic extension of it.
Using this relation we see that the quadratic degree characters can easily be computed
from the knowledge of the inner twists and viceversa.

Let $K_P$ denote the polyquadratic extension of $\Q$ which is the fixed field of the
intersection of all the quadratic degree characters.
It coincides with the kernel of the map $\delta\colon G_\Q\to F^*/F^{*2}$ defined by
putting $\delta(\Frob_p)=\mu_p\pmod{F^{*2}}$ for all primes $p\nmid N$ with $a_p\neq0$
(it is well defined by this condition because the primes $p$ with $a_p\neq0$ are of density one
in the set of all prime numbers for non-CM newforms,
and hence every element of $\Gal(K_P/\Q)$ is a Frobenius for some prime $p$ with $a_p\neq0$).

Let $\sigma_1,\dots,\sigma_r\in G_\Q$ be elements whose restriction to the field $K_P$
is a basis of the Galois group $\Gal(K_P/\Q)$,
and let $\psi_1,\dots,\psi_r\colon G_\Q\to\Z/2\Z$ be (additive versions of)
quadratic degree characters that are a basis of $\Psi$,
and assume that these bases are dual of each other, i.e. such that
$$\psi_i(\sigma_j)=0\quad\text{for}\quad j\neq i\qquad\text{and}\qquad
    \psi_i(\sigma_i)=1,$$
for every $1\leq i,j\leq r$. Let $t_i\in\Q^*$ be rational numbers such that
$$\Q(\sqrt{t_i})=\Qb^{\ker\psi_i}$$
and let
$$\delta_i=\delta(\sigma_i)\in F^*$$
for $i=1,\dots,r$.

We denote by $\gamma_\varepsilon$ the element of the two-torsion
$\Br(\Q)[2]\simeq H^2(G_\Q,\{\pm1\})$
of the Brauer group of the rational numbers that is the cohomology class of the two-cocycle
defined by the formula
$$(\sigma,\tau)\mapsto\sqrt{\varepsilon(\sigma)}\,\sqrt{\varepsilon(\tau)}\,
    {\sqrt{\varepsilon(\sigma\tau)}\,}^{-1}.$$

Then, the Brauer class of $\D=\End^0(B_f)$ is
\begin{equation}\label{brauer-class}
[\D]=\Res_{\Q}^F(\gamma_\varepsilon)\prod(t_i,\delta_i)\in\Br(F)[2],
\end{equation}
and the obstruction to descend the building block over the field $K_P$ is
\begin{equation}\label{obstruction}
\Res_{\Q}^{K_P}(\gamma_\varepsilon)\in\Br(K_P)[2].
\end{equation}

%---------------------------------------------------------------------------------------------------
\section{The package}
%---------------------------------------------------------------------------------------------------

The functions have been implemented in {\tt Magma},
and are based in the packages written by W. Stein for computing with modular symbols,
modular forms, and modular abelian varieties.
In the following we assume that the reader is familiar with {\tt Magma}
and especially with the functions of Stein's packages.

One of the main practical problems in determining the inner twists of a modular form
(which in practice is given by a finite truncation of its Fourier series)
is to check the identity (\ref{definition-psi}) for sufficiently many primes $p$
to guarantee that it holds for all primes.
Since the coefficients for primes smaller than $\frac{k}{12}\psi(N)$
completely determine a cusp form of level $N$ and weight $k$,
with $\psi$ being the function $\psi(n)=n\prod_{p\mid n}(1+1/p)$,
and the twist of a modular form by a Dirichlet character has level dividing the least
common multiple of the level of the form and the square of the conductor of the character,
we see that it is enough to use the bound $\frac k6\psi(N^2)$.
In practice this number becomes too large for computations
even for relatively small levels,
and must be replaced by smaller bounds that should be enough,
but for which there is no proof.
The functions computing the inner twists in W. Stein's package have an optional parameter
to decide which bound to use.
In our implementation we use the large proved bound for levels $N\leq100$
and we use the unproved bound $15+N/2$ for larger levels.
In particular we introduced small modifications in W. Stein's file {\tt inner\_twist.m}.

All the intrinsics we programmed are in a file named {\tt building\_blocks.m};
for using them, just copy that file into the {\tt Magma} packages directory {\tt ModAbVar},
and add the name of the file to {\tt ModAbVar.spec}.
It seems that all these intrinsics will be part of the standard distribution of {\tt Magma} in the future.
The intrinsics implemented are the following:

\begin{itemize}
\item {\tt DegreeMap.} Its input is a space of modular symbols that is new and irreducible,
corresponding to a modular form $f$.
The output consists of a sequence (possibly empty) of pairs $(t_i,\delta_i)$,
with the $t_i$ being nonzero rational numbers only determined modulo squares,
that in fact are normalised to be rational integers that are discriminants of quadratic fields,
and the $\delta_i$ are nonzero numbers of a number field $F$ only determined modulo squares.
These numbers represent the degree map $\delta$ and using them everything else can easily be computed.

\item {\tt BrauerClass.} Its input is a space of modular symbols that is new and irreducible.
The output is the Brauer class of the endomorphism algebra $\D=\End^0(B_f)$,
or equivalently, that of $\End^0(A_f)$, as an element of $\Br(F)[2]$,
computed using the formula (\ref{brauer-class}).
It is given as a sequence with an even number (possibly zero) of the places of the field $F$,
(necessarily finite in weight $2$) in which that quaternion algebra is ramified.
Thus, the variety $B_f$ is a RM-building block (real multiplication)
with $\End^0(B_f)=F$ when this element is trivial,
in which case $A_f\sim B_f^{[E:F]}$ and $\dim B_f=[F:\Q]$;
when this element is nontrivial, then $B_f$ is a QM-building block (quaternionic multiplication)
with $\End^0(B_f)=\D$ a quaternion division algebra over $F$,
in which case $A_f=B_f^{[E:F]/2}$ and $\dim B_f=2[F:\Q]$.

\item {\tt ObstructionDescentBuildingBlock.}
Its input is a space of modular symbols that is new and irreducible.
The output is the obstruction to the existence of a building block over the field $K_P$
that is isogenous to $B_f$.
This obstruction is the element of $\Br(K_P)[2]$ given by (\ref{obstruction})
and is given as the list of the places of the field
for which the local component is nontrivial.
When the list is empty then the field $K_P$ is the smallest possible such field;
otherwise the fields having this property are the extensions of $K_P$ that are splitting fields
for this element of the Brauer group, the smallest possibility being quadratic extensions.

\item {\tt BoundedFSubspace.}
Its input is a Dirichlet character $\varepsilon$, a weight $k$,
and a range of positive integers.
The output is the list of modular symbols spaces of that Nebentypus and weight
corresponding to non-CM newforms, and for which the corresponding field $F$
has degree over Q belonging by the numbers in the given range.

The computation avoids to split the full modular symbols space into newform subspaces,
which is very time-consuming, by the computation of the characteristic polynomials of the
operators $T_p^2/\varepsilon(p)$ for a few primes $p$ not dividing the level,
and the restriction to the modular symbols subspaces that are in the kernel of the
operators obtained from factors of this characteristic polynomial of small degree.
\end{itemize}

%---------------------------------------------------------------------------------------------------
\section{The table}
%---------------------------------------------------------------------------------------------------

Using the functions described in the previous paragraph
we elaborated a table containing information for all non-CM newforms $f$
of weight $2$, level up to $500$, Nebentypus character $\varepsilon$ of not too large order,
and for which the field $F$ has degree over $\Q$ bounded by $4$.
The restrictions introduced in the Nebentypus character are due to the fact that
the Nebentypus modular symbols computations take place in the field generated by
the values of the character, and consequently depend very much on the degree of that field.
So we considered only those characters whose values generate a field of degree up to $12$.

In fact we are quite confident that our table contains all the newforms
in the given level range $N\leq500$, with every Nebentypus character,
and whose corresponding field $F$ has degree up to $4$,
although for checking the nonexistence of such forms with such small $F$
for characters of larger orders,
we replaced the computations of modular symbols over number fields,
that were too slow and required too much memory to be performed,
by computations based in modular symbols over finite fields
and we do not have a full theoretical justification for this reduction process to be correct.

The total number of forms obtained is $5609$.
We give now eight tables containing information for every such form,
corresponding to the eight possible endomorphism algebras structures of $B_f$:
either a totally real field of degree $1,2,3$ or $4$,
or a quaternion division algebra over a totally real field
of dimension $4,8,12$ or $16$ over $\Q$.
For every form we give the following information:
the level $N$, the Nebentypus character $\varepsilon$
(the format is the output of the function {\tt Eltseq}) and its order;
the dimension of the variety $A_f$, equal to the degree of the field $E$ over $\Q$;
information of the field $F$ (when it is not $\Q$) given by a list of coefficients
of a defining polynomial (when it is not quadratic) and its discriminant;
information about the degree map,
given by the numbers $t_i$ and the numbers $N_{F/\Q}(\delta_i)$;
the Brauer class of $A_f$ (if it is nontrivial),
given as a list of integers containing the norms over $\Q$
of the prime ideals of $F$ at which the corresponding quaternion algebra is ramified,
and finally the obstruction to descend $B_f$ to the field $K_P$
(except in the RM-case of even dimension in which we know it is trivial by a theorem of Ribet),
given in the same format but for the field $K_P$.

%\end{document}

{\tiny
\newpage
\section{Tables of newforms corresponding to RM-building blocks}
\subsection{$\Q$-curves}
\setlongtables
\begin{longtable}{|r|c|c||c||c c|}
\hline $N$ & $\varepsilon$ & $\operatorname{ord}(\varepsilon)$ & $\dim(A_f)$ &
$[t_i]$ & $[\delta_i]$ \\
\hline\endhead\hline\endfoot
$11$&$[ 0 ]$&$1$&$1$&$[]$&$[]$\\
$13$&$[ 2 ]$&$6$&$2$&$[ 13 ]$&$[ 3 ]$\\
$14$&$[ 0 ]$&$1$&$1$&$[]$&$[]$\\
$15$&$[ 0, 0 ]$&$1$&$1$&$[]$&$[]$\\
$16$&$[ 0, 1 ]$&$4$&$2$&$[ 8 ]$&$[ 2 ]$\\
$17$&$[ 0 ]$&$1$&$1$&$[]$&$[]$\\
$18$&$[ 2 ]$&$3$&$2$&$[]$&$[]$\\
$19$&$[ 0 ]$&$1$&$1$&$[]$&$[]$\\
\end{longtable}

\newpage
\subsection{RM-surfaces}
\setlongtables
\begin{longtable}{|r|c|c||c||c||c c||c|}
\hline
$N$ & $\varepsilon$ & $\operatorname{ord}(\varepsilon)$ & $\dim(A_f)$ &
$\Disc(F)$ &
$[t_i]$ & $[N_{F/\Q}(\delta_i)]$ & $\gamma_\varepsilon\in\Br(K_P)[2]$ \\
\hline\endhead\hline\endfoot
$17$&$[ 2 ]$&$8$&$4$&$8$&$[ 17 ]$&$[ 2 ]$&$[]$\\
$22$&$[ 2 ]$&$5$&$4$&$5$&$[]$&$[]$&$[]$\\
$23$&$[ 0 ]$&$1$&$2$&$5$&$[]$&$[]$&$[]$\\
$25$&$[ 4 ]$&$5$&$4$&$5$&$[]$&$[]$&$[]$\\
$29$&$[ 0 ]$&$1$&$2$&$8$&$[]$&$[]$&$[]$\\
$31$&$[ 0 ]$&$1$&$2$&$5$&$[]$&$[]$&$[]$\\
$31$&$[ 10 ]$&$3$&$4$&$8$&$[]$&$[]$&$[]$\\
$31$&$[ 6 ]$&$5$&$4$&$5$&$[]$&$[]$&$[]$\\
\end{longtable}

\newpage
\subsection{3-dimensional RM-varieties}
\setlongtables
\begin{longtable}{|r|c|c||c||c|c||c c|}
\hline
$N$ & $\varepsilon$ & $\operatorname{ord}(\varepsilon)$ & $\dim(A_f)$ &
$F$ & $\Disc(F)$ &
$[t_i]$ & $[N_{F/\Q}(\delta_i)]$ \\
\hline\endhead\hline\endfoot
$19$&$[ 2 ]$&$9$&$6$&$[ 0, -3, -1 ]$&$81$&$[]$&$[]$\\
$29$&$[ 4 ]$&$7$&$6$&$[ -5, 6, -1 ]$&$49$&$[]$&$[]$\\
$37$&$[ 4 ]$&$9$&$6$&$[ 0, -3, -1 ]$&$81$&$[]$&$[]$\\
$37$&$[ 4 ]$&$9$&$6$&$[ -6, 9, -1 ]$&$81$&$[]$&$[]$\\
$38$&$[ 2 ]$&$9$&$6$&$[ -6, 9, -1 ]$&$81$&$[]$&$[]$\\
$41$&$[ 0 ]$&$1$&$3$&$[ -1, -3, 1 ]$&$148$&$[]$&$[]$\\
$41$&$[ 10 ]$&$4$&$6$&$[ -1, -5, 3 ]$&$564$&$[ 41 ]$&$[ 2 ]$\\
$43$&$[ 6 ]$&$7$&$6$&$[ -6, 5, -1 ]$&$49$&$[]$&$[]$\\
\end{longtable}

\newpage
\subsection{4-dimensional RM-varieties}
\setlongtables
\begin{longtable}{|r|c|c||c||c|c||c c||c|}
\hline
$N$ & $\varepsilon$ & $\operatorname{ord}(\varepsilon)$ & $\dim(A_f)$ &
$F$ & $\Disc(F)$ &
$[t_i]$ & $[N_{F/\Q}(\delta_i)]$ & $\gamma_\varepsilon\in\Br(K_P)[2]$ \\
\hline\endhead\hline\endfoot
$25$&$[ 2 ]$&$10$&$8$&$[ -11, 36, -31, 1 ]$&$7625$&$[ 5 ]$&$[ 1 ]$&$[]$\\
$32$&$[ 0, 1 ]$&$8$&$8$&$[ 0, -6, -4, 2 ]$&$4352$&$[ 8 ]$&$[ 4 ]$&$[]$\\
$41$&$[ 8 ]$&$5$&$8$&$[ -1, -3, 1, 1 ]$&$725$&$[]$&$[]$&$[]$\\
$41$&$[ 4 ]$&$10$&$8$&$[ -9, 26, -24, 1 ]$&$1125$&$[ 41 ]$&$[ 1 ]$&$[]$\\
$45$&$[ 1, 1 ]$&$12$&$16$&$[ -10, 30, -28, 1 ]$&$54864$&$[ -3, -15 ]$&$[ 1, 4 ]$&$[]$\\
$47$&$[ 0 ]$&$1$&$4$&$[ -1, -5, 5, -1 ]$&$1957$&$[]$&$[]$&$[]$\\
$48$&$[ 0, 1, 0 ]$&$4$&$8$&$[ 0, -6, -4, 2 ]$&$4352$&$[]$&$[]$&$[]$\\
$50$&$[ 4 ]$&$5$&$8$&$[ 1, -9, -4, 16 ]$&$7625$&$[]$&$[]$&$[]$\\
\end{longtable}

\newpage
\section{Tables of newforms corresponding to QM-building blocks}
\subsection{QM-surfaces}
\setlongtables
\begin{longtable}{|r|c|c||c||c c||c||c|}
\hline $N$ & $\varepsilon$ & $\operatorname{ord}(\varepsilon)$ & $\dim(A_f)$ &
$[t_i]$ & $[\delta_i]$ &
$[\D]\in\Br(\Q)[2]$ & $\gamma_\varepsilon\in\Br(K_P)[2]$ \\
\hline\endhead\hline\endfoot
$28$&$[ 1, 1 ]$&$6$&$4$&$[ -7 ]$&$[ 3 ]$&$[ 2, 3 ]$&$[ 2, 2 ]$\\
$35$&$[ 1, 3 ]$&$4$&$4$&$[ -35, -7 ]$&$[ 2, 5 ]$&$[ 2, 5 ]$&$[]$\\
$44$&$[ 1, 5 ]$&$2$&$4$&$[ -11, 44 ]$&$[ 2, 3 ]$&$[ 2, 3 ]$&$[]$\\
$56$&$[ 1, 1, 3 ]$&$2$&$4$&$[ -7 ]$&$[ 6 ]$&$[ 2, 3 ]$&$[ 2, 2 ]$\\
$57$&$[ 1, 9 ]$&$2$&$4$&$[ -19, 57 ]$&$[ 2, 5 ]$&$[ 2, 5 ]$&$[]$\\
$60$&$[ 0, 1, 1 ]$&$4$&$4$&$[ 5, -3 ]$&$[ 2, 5 ]$&$[ 2, 5 ]$&$[]$\\
$63$&$[ 3, 1 ]$&$6$&$4$&$[ -3, -7 ]$&$[ 2, 3 ]$&$[ 2, 3 ]$&$[]$\\
$77$&$[ 3, 5 ]$&$2$&$4$&$[ -11, -7 ]$&$[ 2, 5 ]$&$[ 2, 5 ]$&$[]$\\
\end{longtable}

\newpage
\begin{landscape}
\subsection{4-dimensional QM-varieties}
\setlongtables
\begin{longtable}{|r|c|c||c||c||c c||c||c|}
\hline
$N$ & $\varepsilon$ & $\operatorname{ord}(\varepsilon)$ & $\dim(A_f)$ &
$\Disc(F)$ &
$[t_i]$ & $[N_{F/\Q}(\delta_i)]$ &
$[\D]\in\Br(F)[2]$ & $\gamma_\varepsilon\in\Br(K_P)[2]$ \\
\hline\endhead\hline\endfoot
$36$&$[ 1, 1 ]$&$6$&$8$&$33$&$[ 12, -4 ]$&$[ 1, 3 ]$&$[ 2, 3 ]$&$[]$\\
$40$&$[ 1, 1, 1 ]$&$4$&$8$&$5$&$[ -40, -8 ]$&$[ 4, 20 ]$&$[ 4, 5 ]$&$[]$\\
$76$&$[ 1, 9 ]$&$2$&$8$&$17$&$[ -19, 76 ]$&$[ 8, 13 ]$&$[ 2, 13 ]$&$[]$\\
$88$&$[ 1, 1, 5 ]$&$2$&$8$&$8$&$[ 88, -11 ]$&$[ 49, 56 ]$&$[ 2, 7 ]$&$[]$\\
$99$&$[ 0, 2 ]$&$5$&$8$&$5$&$[ -3 ]$&$[ 11 ]$&$[ 9, 11 ]$&$[]$\\
$108$&$[ 1, 3 ]$&$6$&$8$&$33$&$[ 12, -4 ]$&$[ 1, 3 ]$&$[ 2, 3 ]$&$[]$\\
$120$&$[ 1, 1, 1, 2 ]$&$2$&$16$&$8$&$[ -3, -15 ]$&$[ 8, 56 ]$&$[ 2, 7 ]$&$[]$\\
$135$&$[ 9, 1 ]$&$4$&$8$&$24$&$[ -15, -3 ]$&$[ 1, 10 ]$&$[ 2, 5 ]$&$[]$\\
\end{longtable}

\newpage
\subsection{6-dimensional QM-varieties}
\setlongtables
\begin{longtable}{|r|c|c||c||c|c||c c||c||c|}
\hline
$N$ & $\varepsilon$ & $\operatorname{ord}(\varepsilon)$ & $\dim(A_f)$ &
$F$ & $\Disc(F)$ &
$[t_i]$ & $[N_{F/\Q}(\delta_i)]$ &
$[\D]\in\Br(F)[2]$ & $\gamma_\varepsilon\in\Br(K_P)[2]$ \\
\hline\endhead\hline\endfoot
$48$&$[ 1, 1, 1 ]$&$4$&$12$&$[ -1, -4, 2 ]$&$316$&$[ -24, -3 ]$&$[ 8, 16 ]$&$[ 2, 2 ]$&$[]$\\
$56$&$[ 1, 1, 1 ]$&$6$&$12$&$[ 0, -6, -2 ]$&$756$&$[ 56, -8 ]$&$[ 3, 21 ]$&$[ 3, 7 ]$&$[]$\\
$63$&$[ 1, 3 ]$&$6$&$12$&$[ -7, 10, -3 ]$&$321$&$[ -3, -7 ]$&$[ 3, 27 ]$&$[ 3, 3 ]$&$[]$\\
$77$&$[ 1, 5 ]$&$6$&$12$&$[ -8, 17, -7 ]$&$321$&$[ -7, -11 ]$&$[ 3, 7 ]$&$[ 3, 7 ]$&$[]$\\
$84$&$[ 1, 1, 2 ]$&$6$&$24$&$[ -12, 35, -8 ]$&$2089$&$[ -3, -4 ]$&$[ 8, 64 ]$&$[ 2, 2 ]$&$[]$\\
$91$&$[ 3, 3 ]$&$4$&$12$&$[ -6, 8, -2 ]$&$148$&$[ -91, -7 ]$&$[ 2, 13 ]$&$[ 2, 13 ]$&$[]$\\
$105$&$[ 1, 2, 1 ]$&$6$&$24$&$[ -9, 18, -6 ]$&$1944$&$[ -3, 105, -15 ]$&$[ 6, 12, 28 ]$&$[ 2, 7 ]$&$[]$\\
$129$&$[ 1, 21 ]$&$2$&$12$&$[ -7, 10, -2 ]$&$568$&$[ -3, 129 ]$&$[ 2, 106 ]$&$[ 2, 53 ]$&$[]$\\
\end{longtable}

\newpage
\subsection{8-dimensional QM-varieties}
\setlongtables
\begin{longtable}{|r|c|c||c||c|c||c c||c||c|}
\hline
$N$ & $\varepsilon$ & $\operatorname{ord}(\varepsilon)$ & $\dim(A_f)$ &
$F$ & $\Disc(F)$ &
$[t_i]$ & $[N_{F/\Q}(\delta_i)]$ &
$[\D]\in\Br(F)[2]$ & $\gamma_\varepsilon\in\Br(K_P)[2]$ \\
\hline\endhead\hline\endfoot
$44$&$[ 1, 1 ]$&$10$&$16$&$[ -12, 39, -38, 11 ]$&$7625$&$[ -11, -4 ]$&$[ 1, 11 ]$&$[ 4, 11 ]$&$[]$\\
$76$&$[ 1, 3 ]$&$6$&$16$&$[ -1, -10, 10, -2 ]$&$70164$&$[ -19, 76 ]$&$[ 19, 304 ]$&$[ 2, 19 ]$&$[]$\\
$77$&$[ 3, 1 ]$&$10$&$16$&$[ -8, 14, -7, 1 ]$&$1125$&$[ -11, -7 ]$&$[ 1, 295 ]$&$[ 5, 59 ]$&$[]$\\
$132$&$[ 0, 1, 1 ]$&$10$&$16$&$[ 1, -8, -1, 11 ]$&$5225$&$[ -3, -11 ]$&$[ 121, 3025 ]$&$[ 11, 11 ]$&$[]$\\
$140$&$[ 1, 2, 1 ]$&$6$&$32$&$[ -13, 53, -77, 28 ]$&$90972$&$[ 5, 28, -20 ]$&$[ 28, 324, 1708 ]$&$[ 2, 7 ]$&$[]$\\
$171$&$[ 3, 3 ]$&$6$&$16$&$[ -16, 82, -150, 81 ]$&$248016$&$[ 57, -19 ]$&$[ 36, 81 ]$&$[ 3, 3 ]$&$[]$\\
$195$&$[ 1, 1, 6 ]$&$4$&$32$&$[ -9, 23, -15, 2 ]$&$19796$&$[ 5, -15, -39 ]$&$[ 2, 8, 2300 ]$&$[ 2, 23 ]$&$[]$\\
$225$&$[ 0, 4 ]$&$5$&$16$&$[ -10, 30, -25, 5 ]$&$5125$&$[ -3 ]$&$[ 5 ]$&$[ 5, 9 ]$&$[]$\\
\end{longtable}
\end{landscape}
}
%---------------------------------------------------------------------------------------------------
%---------------------------------------------------------------------------------------------------
\end{document}